\documentclass[12pt]{article}
\usepackage{latexsym,amsfonts,amssymb,amsthm}
\swapnumbers
\newtheoremstyle{notdots}
      {16pt}
      {0pt}
      {\itshape}
      {}
      {\bfseries}
      {}
      {.5em}
      {}

 \theoremstyle{notdots}

\newtheorem{thm}{Theorem}

\newtheorem{prop}[thm]{Proposition}
\newtheorem{lemma}[thm]{Lemma}

\newenvironment{prf}{\noindent\parskip0ex{\textsl{{\bf 
Proof}}}\\}{\hfill\rule{1ex}{1.5ex}\vskip2ex}

\title{A proof of the polycirculant conjecture}
\author{Eric Mwambene, \\Department of Pure and Applied Mathematics,\\
University of the Western Cape,\\Bellville, 7535, South Africa\\
emwambene@uwc.ac.za}
\date{}
\begin{document}
\maketitle
 \begin{abstract}
 This paper presents a solution of the polycirculant conjecture which states that every vertex-transitive graph $G$ has an automorphism that permutes the vertices in cycles of the same length. This is done by identifying vertex-transitive graphs as coset graphs. For a coset graph $H$, an equivalence  relation $\sim$ is defined on the vertices of cosets with classes as double cosets of the stabiliser and any other proper subgroup $A'$ of a transitive group $A$ of $G$. Induced left translations of elements of the subgroup $A'$ are semi-regular since they preserve these double cosets and acts regularly on each of them. The coset graph is equivalent to $G$ by a theorem of Sabidussi.
 
 {\it MSC(2000): 05C25, 20B25}
 \end{abstract}
 \section{Introduction}
This paper sets out to prove  the polycirculant conjecture. The polycirculant conjecture states that every vertex-transitive finite graph is a polycirculant.

A graph $G=(V,E)$ is a set $V$ together with an irreflexive and symmetric  binary relation $E$ defined on $V.$ The two {\it arcs} $(x,y),(y,x)\in E$ are identified into an edge $[x,y]$. All graphs discussed in this paper are loopless, connected and finite.

An automorphism of a graph is a permutation of its vertices which preserves edges. The set of all permutations of a graph $G$ is denoted by Aut $G$ and constitutes a group under composition. A graph $G$ is vertex-transitive if given any two vertices $x,y\in V(G)$ there is an automorphism which maps $x$ to $y.$ Let $r\ge 1$ and $s\ge 2$ be integers. An automorphism $\alpha$ of a graph is said to ($r,s$){\it-semiregular} if  it has $r$ orbits of length $s$ and {\it semiregular} if it is ($r,s$)-semiregular for some $r$ and $s$. A graph $G$ of order $n$ is called a {\it circulant} if it has a (1,n)-semiregular automorphism. 

It is easy to see that a graph is a circulant if and only if it is Cayley. A  natural relaxation of the concept of a graph being a circulant is the concept of it being polycirculant. In 1981, Dragan Maru\v si\v c \cite{marusic} asked if every vertex-transitive graph is a polycirculant. The question was also independently asked by D. Jordan \cite{jordan}. 

The main thrust of this note is the following result.
\begin{thm}
Every vertex-transitive graph is a polycirculant.
\end{thm}

To prove the  theorem, we first represent vertex-transitive graphs as coset graphs. It is a  classical Sabidussi Representation Theorem\cite{sab} that  vertex-transitive graphs are coset graphs.

For convenience and to fix notation, we will first present the Sabidussi Representation Theorem and then we give the proof of our main result.
\section{Main Results}
\subsection{Sabidussi's Representation Theorem}
In this section we give a presentation of Sabidussi Representation Theorem that is amenable to proving  the polycirculant theorem.

Let $G$ be a vertex-transitive graph and $A$ a subgroup of Aut $G$ which acts transitively on $V(G)$.

Fix $u\in V(G)$ (base point) and consider the stabilizer of $u$ in $A$:
$$A_u:=\{\alpha\in A:\alpha(u)=u\}.$$ Let $A/A_u$ be the left cosets of $A_u$ in $A$ and $B:=\{\alpha\in A:[u,\alpha(u)]\in E(G)\}$. 

Note that $B$ is a Cayley set in $A,$ i.e, the identity is not in $B$ and if $b\in B$ then 
$b^{-1}\in B.$

We define a graph $H$ as follows:
$$V(H):=A/A_u,$$ $$[\alpha A_u,\beta A_u]\in E(H)\iff {\rm ~there ~exists~}\upsilon\in 
\beta A_u {\rm~ such~ that ~} \upsilon \in \alpha A_u\sigma, \sigma\in B.$$
\begin{lemma}\label{lemma}
With the notation as above, $H$ is isomorphic to $G.$
\end{lemma}

The constructed graph $H$ is called a {\it coset graph.}
In view of Lemma \ref{lemma} we have that 
\begin{thm}
{\rm (Sabidussi \cite{sab})} Every vertex-transitive graph is a coset graph.
\end{thm}

\subsection{Proof of the theorem}
Here the notation is as in Section 2.1.
Our proof of the theorem utilizes the relationship between Cayley graphs and vertex-transitive graphs. The most used characterization of Cayley graph is that given by Sabidussi.  It reads as follows:
\begin{prop}{\rm (Sabidussi~\cite{sab})}
A graph $G$ is Cayley if and only if {\rm Aut} $G$ contains a subgroup $A$ which acts regularly on $V(G)$. 
\end{prop}
However, for our purpose we present yet another characterization which identifies and uses the fact that Cayley graphs should also be considered as coset graphs.
\begin{prop}\label{prop4}
A graph $G$ is Cayley if and only if there is a transitive subgroup $A$ of {\rm Aut}~$G$ such that $A_u$ (the stabilizer of $u$ in $A$) is normal in $A$ for any $u\in V(G)$.
\end{prop}
\begin{prf}
By Proposition 4, for a given Cayley graph $G,$ {\rm Aut}~$G$ contains a subgroup $A$ which acts regularly on $V(G)$. Moreover, $A_u$ contains only the identity element and therefore is trivially normal in $A$.\vspace{2ex}

Suppose that $G$ has it that Aut $G$ contains a transitive subgroup $A$ such that $A_u$ is normal in $A.$ Then the quotient group $A/A_u$ together with the set 
$$C:=\{\alpha A_u\in A/A_u:\alpha\in A, [u,\alpha u]\in E(G)\}$$ forms a Cayley graph Cay($A/A_u,C$)  isomorphic to $G.$
\end{prf}

To complete the proof of  Theorem 1 we need the following lemma.
\begin{lemma}
Let $G$ be a graph and $u\in V(G)$. Let $A$ be a transitive subgroup of an automorphism group of a graph $G$ and $B:=\{\alpha\in A:[u,\alpha u]\in E(G)\}$. If $A$ contains a proper subgroup $A'$ generated by $B'\subset B$, then it contains a semi-regular element.
\end{lemma}

\begin{prf}
For the vertex-transitive graph $G$, let $H$ be the quotient graph defined as in Section 2.1 by
$$V(H)=A/A_u,$$
$$[\alpha A_u, \alpha' A_u]\in E(H) \iff ~{\rm for ~some~}\tau\in  \alpha A_u,\tau'\in \alpha A_u ~{\rm s.t.}~ \tau\beta=\tau', \beta\in B.$$

We define a relation $\sim$ on $V(H)$ by 
$$\alpha A_u\sim \alpha' A_u\iff \alpha,\alpha'\in A'\tau A_u,$$ i.e. $\alpha A_u$ and $\alpha' A_u$ are related if their union in contained in some double coset $A'\tau A_u$.

It is clear that $\sim$ is an equivalence relation.

We have that for any $\beta\in A'$ the map $\tilde{\lambda_{\beta}}:V(H)\rightarrow V(H)$ given by 
$$\tilde{\lambda_{\beta}}(\alpha A_u)=\beta\alpha A_u$$ is an automorphism of $H$ that preserves equivalence classes of $\sim.$ Moreover, $\Lambda_{A'}:=\{\tilde{\lambda_{\beta}}, \beta\in A'\}$ is regular on each double coset. (We can consider the double coset $A'\tau A_u$ as a set of left cosets in the form $$A'\tau A_u=\{\beta\alpha A_u:\beta\in A'\}$$ and so for each $\beta\in A'$ we have $$\tilde{\lambda_{\beta}}(A'\tau A_u)=\beta A'\tau A_u=A'\tau A_u$$ for each double coset.) Hence the map $\tilde{\lambda_{\beta}}:V(H)\rightarrow V(H)$ is semi-regular.
\end{prf}

{\bf Proof of Theorem 1}\\
Let $G$ be a vertex-transitive graph, $u\in V(G)$ and $B:=\{\alpha\in ~{\rm Aut~} A:[u,\alpha u]\in E(G)\}.$ 

If there is an element $\beta \in B$ such that $<\beta>=$ Aut $G$ 
then the stabilizer of $u$ is normal in Aut $G$, hence by Proposition 5, $G$ is Cayley (Since Aut $G$ is Abelian). We therefore have that $G$ is polycirculant.

Otherwise, $<\beta>, \beta \in B$ is a proper subgroup of Aut $G$ and by Lemma 6, Aut $G$ contains a semi-regular element. (By considering the equivalence relation $\sim$ defining   double cosets in the form $<\beta>\alpha A_u, \alpha\in $ Aut $ G$ as equivalence classes.)\hfill\rule{1ex}{1.5ex}

\end{document}